\documentclass[12pt]{article}
\usepackage{amssymb,amsfonts,amsmath, psfrag,eepic,colordvi}
\topskip  
\numberwithin{equation}{section}
\newtheorem{theorem}{Theorem}[section]

\newtheorem{corollary}[theorem]{Corollary}
\newtheorem{definition}[theorem]{Definition}

\newtheorem{lemma}[theorem]{Lemma}

\newtheorem{example}[theorem]{Example}

\def\pf{\noindent{\it Proof.}}
\def\qed{\nopagebreak\hfill{\rule{4pt}{7pt}}
\medbreak}

\newenvironment{kst}
{\setlength{\leftmargini}{1.1\parindent}
\begin{itemize}
\setlength{\itemsep}{-2.0mm}} {\end{itemize}}

\makeatletter
\long\def\@makecaption#1#2{%
   \vskip 10\p@
   \setbox\@tempboxa\hbox{{#1}.\ \ #2}%
   \ifdim \wd\@tempboxa >\hsize
       {#1}\ \ #2\par
   \else
       \hbox to\hsize{\hfil\box\@tempboxa\hfil}%
   \fi}
\makeatother

\makeatletter
\renewcommand{\@seccntformat}[1]{{\csname the#1\endcsname}{\normalsize .}\hspace{.5em}}
\makeatother

\begin{document}
\parskip 6pt

\pagenumbering{arabic}
\def\sof{\hfill\rule{2mm}{2mm}}
\def\ls{\leq}
\def\gs{\geq}
\def\SS{\mathcal S}
\def\qq{{\bold q}}
\def\MM{\mathcal M}
\def\TT{\mathcal T}
\def\EE{\mathcal E}
\def\lsp{\mbox{lsp}}
\def\rsp{\mbox{rsp}}
\def\pf{\noindent {\it Proof.} }
\def\mp{\mbox{pyramid}}
\def\mb{\mbox{block}}
\def\mc{\mbox{cross}}
\def\qed{\hfill \rule{4pt}{7pt}}
\def\block{\hfill \rule{5pt}{5pt}}

\begin{center}
{\Large\bf Matchings Avoiding Partial Patterns} \vskip 6mm

William Y. C. Chen$^1$, Toufik Mansour$^{2,1}$, Sherry H. F.
Yan$^3$

\vskip 3mm $^{1,3}$Center for Combinatorics, LPMC, Nankai
University\\
  Tianjin 300071, P.R. China

$^2$Department of Mathematics, University of Haifa, 31905 Haifa,
Israel.

{\tt $^1$chen@nankai.edu.cn, $^2$toufik@math.haifa.ac.il,
$^3$huifangyan@eyou.com}
\end{center}
\vskip 2mm
\begin{center}
 ABSTRACT
\end{center}
We show that  matchings avoiding certain
partial patterns are counted by the $3$-Catalan numbers. We give a
characterization of $12312$-avoiding matchings in terms of
restrictions on the corresponding oscillating tableaux. We also
find a bijection between Schr\"oder paths without peaks at level
one and matchings avoiding both patterns $12312$ and $121323$.
Such objects are counted by the super-Catalan numbers or the
little Schr\"{o}der numbers. A refinement of the super-Catalan
numbers is obtained by fixing the number of crossings in the
matchings. In the sense of Wilf-equivalence, we find that the
patterns 12132, 12123, 12321, 12231, 12213  are equivalent to
$12312$.
\medskip

\noindent {\sc Key words}: Generating function, generating tree,
matching, ternary tree, super-Catalan number, oscillating tableau.

\noindent {\sc AMS Mathematical Subject Classifications}: 05A05,
05C30.


\section{\normalsize Introduction}

A {\it matching} on a set $[2n]=\{1,2,\ldots,2n\}$ is a partition
of $[2n]$ in which every block contains exactly two elements, or
equivalently a graph on $[2n]$ in which every vertex has degree
one. There are many ways to represent a matching. It can be
displayed by drawing the $2n$ points on a horizontal line in the
increasing order. This  is called the {\em linear representation}
of a matching \cite{c3}. An edge $(i,j)$ is drawn as an arc
between the nodes $i$ and $j$ above the horizontal line, where the
vertices $i$ and $j$ are called the initial point and the end
point, respectively. An edge $e=(i,j)$ is always written in such a
way that $i<j$. Let $e=(i,j)$ and $e'=(i',j')$ be two edges of a
matching $M$, we say that $e$ {\it crosses} $e'$ if they intersect
with each other, in other words, if $i<i'<j<j'$. In this case, the
pair of edges $e$ and $e'$ is called a {\em crossing} of the
matching. Otherwise, $e$ and $e'$ are said to be {\it
noncrossing}. The set of matchings on $[2n]$ is denoted by
$\MM_n$. Note that $|\MM_n|=(2n-1)!!=1\cdot3\cdot5\cdots(2n-1)$.

\begin{figure}[h,t]
\begin{center}
\begin{picture}(15,3)
\setlength{\unitlength}{2mm} \linethickness{0.4pt}
\qbezier(0,0)(3,5)(6,0)\qbezier(2,0)(5,5)(8,0)\qbezier(4,0)(7,5)(10,0)
\put(4, -1.5){\small$3$}\put(1.8, -1.5){\small$2$}\put(-0.2,
-1.5){\small$1$}\put(5.8, -1.5){\small$4$}\put(7.7,
-1.5){\small$5$}\put(9.6, -1.5){\small$6$}
\end{picture}
\end{center}
\caption{The matching $123123$.} \label{fig1}
\end{figure}
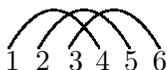

In this paper, we also use the representation of a matching $M$ of
$n$ edges by a sequence of length $2n$ on the set $\{1, 2, \ldots,
n\}$ such that each element $i$ $(1\leq i \leq n)$ appears exactly
twice, and the first occurrence of the element $i$ precedes that
of  $j$ if $i<j$.  Such a representation is called the {\em
Davenport-Schinzel sequence} \cite{4,13} or the {\it canonical
sequential form} \cite{10}. In fact, the canonical sequential
representation of a matching is the sequence obtained from its
linear representation by labeling the endpoints of each arc in the
order of the appearance of its initial point such that the
endpoints of each arc have the same label. For example, the
matching in Figure~\ref{fig1} can be represented by $123123$.

Let $\pi=\pi_1\pi_2 \ldots \pi_k$ and
$\tau=\tau_1\tau_2\ldots\tau_k$ be two sequences. If for any
$1\leq i,j\leq k$ we have $\pi_i<\pi_j$ if and only if
$\tau_i<\tau_j$, then we say $\pi$ and $\tau$ are {\em
order-isomorphic}. The matching $\pi$ contains an occurrence of
$\tau$ if there is a subsequence in the canonical sequential form
of $\pi$ which is order-isomorphic to $\tau$. In such a context
$\tau$ is usually called a \emph{pattern}. When a pattern forms a
representation of a small matching, we say that it is complete;
otherwise, we say that it is partial.  In this paper we are mainly
concerned with the partial pattern $12312$. We say that $\pi$ {\em
avoids} $\tau$, or $\tau$-{\it avoiding}, if there is no
occurrence of the pattern $\tau$ in the matching $\pi$. The set of
all $\tau$-avoiding matchings on $[2n]$ is denoted $\MM_n(\tau)$.
Denote by $M_n(\tau_1, \tau_2, \ldots, \tau_k)$ the set of
matchings on $[2n]$ which avoid the patterns $\tau_1, \tau_{2},
\ldots, \tau_k$. Pattern avoiding matchings  have been studied by
de M\'edicis and Viennot \cite{14}, de Sainte-Catherine \cite{7},
Gessel and Viennot \cite{5}, Gouyou-Beauchamps \cite{8,9}, Stein
\cite{17}, Touchard \cite{18}, and recently by Klazar
\cite{10,11,12}, Chen, Deng, Du, Stanley and Yan \cite{c}.

The {\em $k$-Catalan numbers}, or generalized Catalan numbers are
defined by
         $$C_{n,k}={1\over {(k-1)n+1}}{kn\choose n}$$
for $n\geq 1$ (see \cite{hp}). For $k=2$, the $2$-Catalan numbers
are the usual Catalan numbers.

In this paper we show that  $12312$-avoiding matchings on $[2n]$
are  counted by the $3$-Catalan number, namely,
       $$|\MM_n(12312)|=\frac{1}{2n+1}\binom{3n}{n}.$$
We note that the following objects are also counted by the
$3$-Catalan numbers:
\begin{kst}
\item complete ternary trees with $n$ internal nodes, or $3n$
edges \cite{pp},

\item even trees with $2n$ edges \cite{c1, df}, \item  noncrossing
trees with $n$ edges \cite{dn, pp}, \item the set of lattice paths
from $(0,0)$ to $(2n, n)$ using steps $E=(1,0)$ and $N=(0,1)$ and
never lying above the line $y=x/2$ \cite{hp},

\item dissections of a convex $2n+2$-gon into $n$ quadrilaterals
by drawing $n-1$ diagonals, no two of which intersect in its
interior \cite{hp},

\item two line arrays $\binom{\alpha}{\beta}$, where
$\alpha=\{a_1, a_2, \ldots, a_n\}$ and $\beta=\{b_1, b_2, \ldots,
b_n\}$ such that  $1=b_1=a_1\leq b_2\leq a_2 \ldots \leq b_n\leq
a_n$ and $a_i\leq i$ \cite{ca}.
\end{kst}
The relations between  ternary trees, even trees, and noncrossing
trees have been studied by Chen \cite{c1}, Feretic and Svrtan
\cite{fs}, Noy \cite{n}, and Panholzer and Prodinger \cite{pp}.
Stanley discussed several of these families in \cite[Problems
$5.45-5.47$]{st}.

By using generating functions, we derive a formula for the number
of matchings in $\MM_n(12312)$ having exactly $m$ crossings. We
also show that the cardinality of $\MM_{n-1}(12312, 121323)$ is
the $n$-th super-Catalan number or the little Schr\"{o}der number
for $n\geq 1$ (see~\cite[Sequence A001003]{Seq}). By considering
the number of matchings in $\MM_{n-1}(12312, 121323)$ having
exactly $m$ crossings we obtain a closed expression for a
refinement of the super-Catalan numbers. The $n$-th super-Catalan
number also counts the number of Schr\"oder paths of semilength
$n-1$ (i.e. lattice paths from $(0,0)$ to $(2n-2,0)$, with steps
$H=(2,0)$, $U=(1,1)$, and $D=(1,-1)$ and not going below the
$x$-axis) with no peaks at level one, as well as certain Dyck
paths (see~\cite[Sequence A001003]{Seq} and references therein).
We find a bijection between Schr\"oder paths of semilength $n$
without peaks at level one and matchings on $[2n]$ avoiding both
patterns $12312$ and $121323$.

Following the approach of Chen, Deng, Du, Stanley and Yan
\cite{c}, we use oscillating tableaux to study $12312$-avoiding
matchings. The notion of {\em oscillating tableaux} is introduced
by Sundaram~\cite{Sun86,su} in the study of the representations of
the symplectic group (see also \cite{De1,Roby}). These tableaux
play an important role in Berele's decomposition
formula~\cite{Be1} for powers of defining representations of the
complex symplectic groups. In fact, an oscillating tableau is a
sequence of Young diagrams (or partitions) starting and ending
with the empty diagram $\lambda: \emptyset=\lambda^0, \lambda^1,
\ldots \lambda^{k-1},\lambda^{k}=\emptyset$ such that the diagram
$\lambda^i$ is obtained from $\lambda^{i-1}$ by either adding one
square or removing one square. An oscillating tableau can be
equivalently formulated as a sequence of standard Young tableaux
(often abbreviated as SYT). The number $k$ in the above definition
is called the length of the oscillating tableau $\lambda$.

It has been shown by Stanely~\cite{st} that oscillating tableaux
of length $2n$ are in one-to-one correspondence with matchings on
$[2n]$. In this paper we apply this bijection to $12312$-avoiding
matchings and obtain the corresponding oscillating tableaux and
closed lattice walks. We further provide a one-to-one
correspondence between the set of closed lattice walks and the set
of lattice paths from $(0,0)$ to $(2n, n)$ using steps $E=(1,0)$
and $N=(0,1)$ without crossing the line $y=x/2$, see \cite{GJ}.
From this perspective, we see that $\MM_n(12312)$ is counted by
the $3$-Catalan numbers.

In addition to the pattern $12312$, we find other  patterns that
are equivalent to $12312$ in the sense of Wilf-equivalence. To be
more specific, we show that for any pattern $\tau \in \{12312,
12132, 12123, 12321, 12231, 12213\}$, we have
                     $|M_{n}(\tau)|=C_{n,3}$.
We use the technique of generating trees to reach this conclusion.
A generating tree is a rooted tree in which each node is
associated with a label, and the labels of the children of a node
are determined by certain succession rules. The idea of generating
trees was introduced by Chung, Graham, Hoggat, JR. and M. Kleiman
\cite{3} for the study of Baxter permutations and was further
applied to the study of pattern avoidance by Stankova and West
\cite{15,16,19,20}.  Barcucci et al. \cite{2} developed the ECO
method: a methodology for enumeration of combinatorial objects,
which is based on the technique of generating trees.

\section{\normalsize Matchings and ternary trees}

In this section, we use the sequence representation of a matching
as described in the introduction. Our goal is to show that
$\MM_{n}(12312)$ is counted by the $3$-Catalan number. The first
approach is to give a recursive construction of the set
$\MM_{n}(12312)$. Intuitively, for a matching $\theta=a_1a_2\cdots
a_{2n}$ in $\MM_{n}(12312)$, we may obtain a matching in
$\MM_{n-1}(12312)$ by removing an edge. Then we need to keep track
of all possible ways to recover a matching in $\MM_{n}(12312)$
from a smaller matching. In the recursive generation of matchings
with $n$ edges, one is often concerned with the edge whose initial
and end points have the label $n$ in the canonical sequential
form. However, for the purpose of this paper, we use the edge that
is associated with the last node $2n$.  We denote by $E_\theta$
the edge $(j, 2n)$ that is associated with the last node $2n$. In
general, we use the notation $E_i$ to denote the edge with end
point $i$. In this sense, $E_\theta=E_{2n}$.

Let $E_\theta=(j, 2n)$ be the edge of $\theta$ associated with the
last node $2n$. Clearly, if $\theta$ is $12312$ avoiding, then the
matching $\theta'$ obtained from $\theta$ by removing the edge
$E_\theta$ is also $12312$-avoiding. Thus, the question becomes
how to identify the possibilities of the position $j$ in the
matching $\theta'$ for which one can add the edge $(j,2n)$ to form
a $12312$-avoiding matching.

We need to introduce the notion of the critical crossing of a
matching $\theta$. Let $F_\theta$ be the edge with the rightmost
end point that intersects with the edge $E_\theta$. We call
$F_\theta$ the critical edge of $\theta$.  If $\theta$ is
$12312$-avoiding, then the subgraph induced by the nodes between
the end points of $F_\theta$ and $E_\theta$ is a $12312$-avoiding
matching. If there does not exist any edge that intersects with
$E_\theta$, then the subgraph induced by the nodes between the
initial point and the end point of $E_\theta$ is a
$12312$-avoiding matching.

Let us now consider the case when there exists a critical
crossing. We have the following lemma on the structure of
$12312$-avoiding matchings, which is straightforward to verify.

\begin{lemma}\label{lem1}
Let $\theta$ be a $12312$-avoiding matching that has a critical
crossing. Let $E_\theta=(j,2n)$ and $F_\theta=(x,y)$. Then all the
nodes between $j$ and $y$ are end points, and the edges associated
with these nodes do not cross each other and have their initial
nodes between $x$ and $j$.
\end{lemma}
 Suppose that $\theta$ has a critical crossing. Let $i$
be initial point of the edge $E_{j+1}$. As a consequence of the
above lemma, we see that the subgraph induced by the nodes between
$i$ and $j$ forms a $12312$-avoiding matching. It remains to
consider the structure of the edges associated with the nodes
before the node $i$. We need the following observation.
\begin{lemma}\label{lem2}
Let $\theta$ be a $12312$-avoiding matching that has a critical
crossing. Let $E_\theta=(j,2n)$ and $F_\theta=E_{j+m}$. Then
$E_\theta$ is the only edge that intersects any two edges
$E_{j+r}$ and $E_{j+s}$ for $1\leq r<s\leq m$.
\end{lemma}

The above lemmas are sufficient to demonstrate the recursive
structure of $12312$-avoiding matchings. We need to decompose the
 matching on  $\{1, 2, \ldots, j\}$ into segments for the
construction of smaller $12312$-avoiding matchings. The first step
is to find an edge that intersects with $F_\theta$ with the
rightmost end point $v$. If such edge does not exist then we get a
$12312$-avoiding matching induced by the nodes from $1$ to the
node before the initial point of $F_\theta$. Otherwise, we have an
edge with rightmost end point $v$ that intersects with $F_\theta$.

We claim that the nodes from $1$ to $v$, altogether with the node
$j+m$ (the end point of $F_\theta$), form a $12312$-avoiding
matching. This can be seen from the fact that all the nodes
between the initial point of $F_\theta$ and $v$ are end points
since $\theta$ is $12312$-avoiding.

Now we ready to describe the recursive structure of
$12312$-avoiding matchings. If a $12312$-avoiding matching does
not have a critical crossing, then it consists of two smaller
$12312$-avoiding matchings as illustrated by Figure~\ref{rec11}.
\begin{figure}[h,t]
\begin{center}
\begin{picture}(18,6)
\setlength{\unitlength}{2mm} \linethickness{0.4pt}
\qbezier(2,0)(5,4)(8,0)\put(2,0){\circle*{0.4}}\put(8,0){\circle*{0.4}}
\put(0,0){\footnotesize{\mbox{$\blacksquare$}}}
\put(4.5,0){\footnotesize{\mbox{$\blacksquare$}}}
\end{picture}
\end{center}
\caption{} \label{rec11}
\end{figure}

The nontrivial part is the recursive structure of $12312$-avoiding
matchings that have critical crossings. Let $\theta$ be a matching
on $[2n]$ having the two edges $E_\theta=(j,2n)$ and
$F_\theta=(i,j+m)$. Let us consider the following two cases: (1)
There exists an edge $E_y$ with the rightmost end point  crossing
the edge $F_\theta$. (2) There does not exist such an edge.

For the above two cases, we see that $\theta$ can always be
decomposed into smaller $12312$-avoiding matchings.   In the first
case, we obtain three smaller matchings given below.
\begin{kst}
\item[(1A)] The induced subgraph of $\theta$ on the nodes $1,2,\ldots,y,j+m$.\\[-7pt]
\item[(1B)] The induced subgraph of $\theta$ on the nodes $y+1,\ldots,j+m-1,2n$.\\[-7pt]
\item[(1C)] The induced subgraph of $\theta$ on the nodes
$j+m+1,\ldots,2n-1$.
\end{kst}
In the second case, we obtain three smaller matchings:
\begin{kst}
\item[(2A)] The induced subgraph of $\theta$ on the nodes $1,2,\ldots, i, j+m$.\\[-7pt]
\item[(2B)] The induced subgraph of $\theta$ on the nodes $i+1,\ldots,j+m-1,2n$.\\[-7pt]
\item[(2C)] The induced subgraph of $\theta$ on the nodes
$j+m+1,\ldots,2n-1$.
\end{kst}

To find the number of $12312$-avoiding matchings on $[2n]$ with a
given number of crossings, we need a refinement of the above
structure of $12312$-avoiding matchings. We see that the matching
$\theta$ avoids $12312$ if and only if the smaller structures
$\theta_1$ and $\alpha'$ are $12312$-avoiding matchings, where
$\theta_1$ is the corresponding induced subgraph in Case (1A) or
(2A) (in Case (1) we define $\theta_1$ as the empty matching), and
$\alpha'$ is the induced subgraph obtained from $\theta$ by
deleting the subgraph $\theta_1$. The matching $\theta_1$ has the
critical edge of $\theta$, namely $E_{j+m}$.
\begin{figure}[h,t]
\begin{center}
\begin{picture}(100,15)
\setlength{\unitlength}{2mm} \linethickness{0.4pt}
\qbezier(0,0)(20,15)(40,0)\put(0,0){\circle*{0.4}}\put(40,0){\circle*{0.4}}
\qbezier(8,0)(20,8)(32,0)\put(8,0){\circle*{0.4}}\put(32,0){\circle*{0.4}}
\qbezier(12,0)(20,4)(28,0)\put(12,0){\circle*{0.4}}\put(28,0){\circle*{0.4}}
\put(4,0){\small{$\ldots$}}\put(34,0){\small{$\ldots$}}
\qbezier(24,0)(39,15)(54,0)\put(24,0){\circle*{0.4}}\put(54,0){\circle*{0.4}}
\put(-0.2,-0.9){\footnotesize{\mbox{$\blacksquare$}}}\put(-0.2,-2.2){\footnotesize{$\theta_1$}}
\put(7.8,-0.9){\footnotesize{\mbox{$\blacksquare$}}}\put(7.6,-2.2){\footnotesize{$\theta_{m-1}$}}
\put(11.8,-0.9){\footnotesize{\mbox{$\blacksquare$}}}\put(11.9,-2.2){\footnotesize{$\theta_m$}}
\put(18,-0.7){\footnotesize{\mbox{$\blacksquare$}}}\put(19,-0.7){\footnotesize{\mbox{$\blacksquare$}}}\put(18.7,-1.8){\footnotesize{$\alpha$}}
\put(45,-0.7){\footnotesize{\mbox{$\blacksquare$}}}\put(46,-0.7){\footnotesize{\mbox{$\blacksquare$}}}\put(45.7,-1.9){\footnotesize{$\beta$}}
\put(26.5,-1.5){\footnotesize{$j+1$}}\put(30.6,-1.5){\footnotesize{$j+2$}}\put(38.5,-1.5){\footnotesize{$j+m$}}
\put(23.5,-1.5){\footnotesize{$j$}}\put(53,-1.5){\footnotesize{$2n$}}
\end{picture}
\end{center}
\caption{} \label{structure1}
\end{figure}

Now, let us consider $\alpha'$ which is a smaller $12312$-avoiding
matching. If there exists a critical crossing, say $E_{j+m-1}$,
for the matching $\alpha'$, then the matching $\alpha'$ can be
decomposed into two smaller $12312$-avoiding matchings $\theta_2$
and $\alpha''$. Repeating the above procedure $m$ times we see
that the matching $\theta$ can be decomposed into $m+2$ smaller
$12312$-avoiding matchings where the first $m$ smaller matchings
  have critical edges $E_{j+m},\ldots,E_{j+1}$. We call
these $m$ edges quasi-critical edges of $\theta$
(see~Figure~\ref{structure1}).
\begin{lemma}\label{blem1}
Let $\theta$ be any matching on $[2n]$ with $E_\theta=(j,2n)$.
Then, there exist $m$ quasi-critical edges
$E_{j+1},E_{j+2},\ldots,E_{j+m}$ of $\theta$ such that $\theta$
can be decomposed into $m+2$ smaller $12312$-avoiding matchings
$\theta_1,\ldots,\theta_m,\alpha,\beta$ such that
\begin{kst}
\item[{\rm(1)}] $\theta_s$ is the induced subgraph on the nodes
$v_{s-1}+1,v_{s-1}+2,\ldots,v_s,j+m+1-s$,

\item[{\rm(2)}] $\alpha$ is the induced subgraph on the nodes
$v_m+1,v_m+2,\ldots,j-1$,

\item[{\rm(3)}] $\beta$ is the induced subgraph on the nodes
$j+m+1,j+m+2,\ldots,2n-1$,
\end{kst}
where $v_0=0$ and $v_r$, $r=1,2,\ldots,m$, is the rightmost end
point of an edge that crosses the edge $E_{j+m+1-r}$. If such an
edge does not exist, we define $v_r$ as the initial point of
$E_{j+m+1-r}$.
\end{lemma}

Moreover, as a corollary of Lemma~\ref{blem1} we find a formula
for the number of $12312$-avoiding matchings on $[2n]$ with
exactly $m$ crossings.

\begin{theorem}\label{mth1}
The number of $12312$-avoiding matchings  on $[2n]$ with exactly
$m$ crossings is given by
$$\sum_{i=n}^{2n-1}\frac{(-1)^{n+m+i}}{i}\binom{i}{n}\binom{3n}{i+1+n}\binom{i-n}{m}.$$
\end{theorem}
\pf Let
\[ G(x,y)=\sum_{n\geq 0}\; \sum_{\theta\in
M_{n}(12312)}x^ny^{c(\theta)},\] where $c(\theta)$ is the number
of crossings of $\theta$. Let
\[B(x,y)=\sum_{n\geq
1}\sum_{\theta}x^ny^{c(\theta)},\]
 where the second summation ranges
over  matchings $\theta_s$ as in the first case of
Lemma~\ref{blem1}. It follows from Lemma~\ref{blem1}  that the
ordinary generating function for the number of $12312$-avoiding
matchings with  exactly $m$ quasi-critical  edges
$E_{j+1},\ldots,E_{j+m}$ is given by $xy^mG^2(x,y)B^{m}(x,y)$.
Summing over all the possibilities for $m\geq0$ we arrive at
\begin{equation}\label{e1}
G(x, y)=1+{{xG^2(x, y)}\over{1-yB(x, y)}}.
\end{equation}
Applying Lemma~\ref{blem1} for matchings of the form
$\theta_j$, it follows that that the ordinary
generating function for the number of $12312$-avoiding matchings
$\theta_j$ with exactly $k$ quasi-critical edges  is given by
$xy^kG(x,y)B^{k}(x,y)$. Therefore, summing over all the
possibilities for $k\geq0$ we get
\begin{equation} \label{e2}
         B(x, y)={xG(x,y)\over {1-yB(x,y)}}.
\end{equation}
Combining (\ref{e1}) and (\ref{e2}) we get
\begin{equation}\label{e3}
 B(x,y)={{G(x,y)-1}\over G(x,y)}.
 \end{equation}
It follows from (\ref{e1}) and (\ref{e3}) that $G(x,y)$ satisfies
the following recurrence relation
\begin{equation} \label{e4}
xG(x,y)^3+G(x,y)-G(x,y)^2+y(G(x,y)-1)^2=0.
\end{equation}
Substituting $xy$ by $x$ and $y+1$ by $y$ in above recurrence, we
get
\begin{equation} \label{e5}
G(xy,y+1)=1+y(xG^3(xy,y+1)+(G(xy,y+1)-1)^2).
\end{equation}
Using the Lagrange inversion formula we obtain
$$G(xy,y+1)=1+\sum\nolimits_{i\geq1}\frac{1}{i}\sum\nolimits_{j=0}^i\binom{i}{j}\binom{3j}{i+1+j}x^jy^{i},$$
which implies that
\begin{equation}\label{gxy}
G(x,y)=1+\sum\nolimits_{i\geq1}\frac{1}{i}
\sum\nolimits_{j=0}^i\binom{i}{j}\binom{3j}{i+1+j}x^j(y-1)^{i-j}.
\end{equation}
Then $[x^ny^m]G(x,y)$ gives the number of $12312$-avoiding
matchings on $[2n]$ with exactly $m$ crossings. \qed

Setting $y=1$ in (\ref{gxy}), we obtain the following conclusion.

\begin{theorem}\label{theo1}
The number of $12312$-avoiding matchings on $[2n]$ equals the $3$-Catalan
number $C_{n,3}$.
\end{theorem}

In fact, we may use the above recursive structure of
$12312$-avoiding matchings to construct a bijection between
$12312$-avoiding matchings on $[2n]$ and ternary trees with $n$
internal nodes.  Instead, we will construct a bijection between
$12312$-avoiding matchings and oscillating tableaux which are in
one-to-one correspondence with lattice paths counted by $C_{n,3}$.
Let us consider the case $m=0$, namely, the number of
$12312$-avoiding matchings without any crossings. It is clear that
any noncrossing matching automatically avoids the pattern $12312$.
Therefore, the above formula reduces to the Catalan number when
$m=0$.

\begin{corollary}
For all $n\geq 1$, we have the identity
$$\sum_{i=n}^{2n-1}\frac{(-1)^{n+i}}{i}\binom{i}{n}\binom{3n}{i+1+n}={1\over
{n+1}}{2n\choose n}.$$
\end{corollary}

\section{\normalsize $\MM_n(12312, 121323)$ and Schr\"oder paths }

In this section, we are concerned with the matchings avoiding both
patterns $12312$ and $121323$. We need a  refinement of
Lemma~\ref{blem1}.

\begin{figure}[h,t]
\begin{center}
\begin{picture}(100,15)
\setlength{\unitlength}{2mm} \linethickness{0.4pt}
\qbezier(0,0)(20,15)(40,0)\put(0,0){\circle*{0.4}}\put(40,0){\circle*{0.4}}
\qbezier(8,0)(20,8)(32,0)\put(8,0){\circle*{0.4}}\put(32,0){\circle*{0.4}}
\qbezier(12,0)(20,4)(28,0)\put(12,0){\circle*{0.4}}\put(28,0){\circle*{0.4}}
\put(3,0){\small{$\ldots$}}\put(34,0){\small{$\ldots$}}
\qbezier(24,0)(39,15)(54,0)\put(24,0){\circle*{0.4}}\put(54,0){\circle*{0.4}}
\put(-2,0){\footnotesize{\mbox{$\blacksquare$}}}\put(-2,-1.5){\footnotesize{$\theta_1$}}
\put(6.5,0){\footnotesize{\mbox{$\blacksquare$}}}\put(5,-1.5){\footnotesize{$\theta_{m-1}$}}
\put(10.5,0){\footnotesize{\mbox{$\blacksquare$}}}\put(10.5,-1.5){\footnotesize{$\theta_m$}}
\put(18,0){\footnotesize{\mbox{$\blacksquare$}}}\put(19,0){\footnotesize{\mbox{$\blacksquare$}}}\put(17.5,-1.5){\footnotesize{$\theta_{m+1}$}}
\put(45,0){\footnotesize{\mbox{$\blacksquare$}}}\put(46,0){\footnotesize{\mbox{$\blacksquare$}}}\put(45.7,-1.5){\footnotesize{$\beta$}}
\put(26.5,-1.5){\footnotesize{$j+1$}}\put(30.6,-1.5){\footnotesize{$j+2$}}\put(38.5,-1.5){\footnotesize{$j+m$}}
\put(23.5,-1.5){\footnotesize{$j$}}\put(53,-1.5){\footnotesize{$2n$}}
\end{picture}
\end{center}
\caption{} \label{sch1}
\end{figure}

\begin{lemma}\label{prop.1}
Let $\theta$ be a  matching on $[2n]$ with $n\geq 1$. Assume that $\theta$
has $m$ $(m\geq 0)$ quasi-critical edges $E_{j+1},\ldots,E_{j+m}$.
Then $\theta$ avoids both
patterns $12312$ and $121323$ if and only if
$\theta$ can be decomposed into smaller matchings
$\theta_1,\ldots,\theta_{m+1},\beta$ avoiding both patterns
$12312$ and $121323$ such that
\begin{kst}
\item[{\rm{(1)}}] $\theta_s$ is the induced subgraph of $\theta$
on the nodes $v_{s-1}+1,\ldots,v_s-1$,

\item[{\rm{(2)}}] $\beta$ is the induced subgraph of $\theta$ on
the nodes $j+m+1,\ldots,2n-1$,
\end{kst}
where $v_0=0$, $v_{m+1}=j$, and $v_s$ is the initial point of the
edge $E_{j+m+1-s}$.
\end{lemma}

Let
\[ F(x)=\sum_{n\geq 0}f_nx^n\]
be the ordinary generating
function of the number of matchings on $[2n]$ which avoid both
patterns $12312$ and $121323$. From Lemma \ref{prop.1} we have
 the following recurrence relation
            \[ F(x)=1+{xF^2(x)\over {1-xF(x)}}. \]
It follows that
$$F(x)=\frac{1+x-\sqrt{1-6x+x^2}}{4x}
=1+\sum_{n\geq1}\frac{1}{n}\sum_{j=1}^{n}2^{j-1}\binom{n}{j}\binom{n}{j-1}x^n. $$
Now we see that for $n\geq 1$,  $f_{n-1}$ turns out to be
 the $n$-th super-Catalan number which equals the number of
Schr\"oder paths of semilength $n-1$ without peaks at level one.

We proceed to give a bijection  $\phi$  between the set of
Schr\"oder paths of semilength $n$ without peaks at level
 one and the set of  matchings on $[2n]$ which avoid both patterns
$12312$ and $121323$.
Note that any nonempty Schr\"oder path $P$  has the following unique decompostion:
                $$P=HP' \ \mbox{or}\ P=UP'DP'',$$
where $P'$ and $P''$ are possibly empty  Schr\"oder paths. This is
called the \emph{first return decomposition} by Deutsch \cite{De}.

Given a Schr\"oder path $P$ of semilength $n$ without peaks at
level one, if it is empty, then $\phi(P)$ is the empty matching.
Otherwise,  we may decompose it by using the first return
decomposition. We may use this decomposition recursively to
get  the matching $\phi(P)$ on $[2n]$
avoiding both patterns $12312$ and $121323$.
\begin{kst}
\item [Case 1.] If $P=HP'$,  we have the structure as shown in
Figure \ref{case1}:
\begin{figure}[h,t]
\begin{center}
\begin{picture}(10,5)
\setlength{\unitlength}{4mm} \linethickness{0.4pt}
\put(-4, -0.2){\small$\phi(P)$}\put(-1.5,
-0.1){$=$}\qbezier(0,0)(1.5,2)(3,0)\put(1.5,0){$\block$}\put(1,
-1){\small{$\phi(P')$}}
\end{picture}
\end{center}
\caption{Case 1.}\label{case1}
\end{figure}

\item[Case 2.] If $P=UP'DP''$ and $P'=P_1UDP_2UD\ldots
P_kUDP_{k+1}$, where for any $1\leq i\leq k+1$, $P_i$ is a
Schr\"oder path without peaks at level one,  then we have the
structure as shown in Figure \ref{case2}:
\begin{figure}[h,t]
\begin{center}
\begin{picture}(50,10)
\setlength{\unitlength}{4mm} \linethickness{0.4pt}
\put(-4, -0.2){\small$\phi(P)$}\put(-1.5, -0.1){$=$}
\put(0,0){$\block$}\put(-0.8, -1){\small$\phi(P_1)$}
\qbezier(1,0)(6,3.5)(11,0)\put(2,0){$\block$}\put(1.5,
-1){\small$\phi(P_2)$}\qbezier(3,0)(6,2)(9,0)\put(3.5,
0){$\ldots$}\put(6,0){$\block$}\put(5.5,
-1){\small$\phi(P_{k+1})$} \qbezier(7,0)(11,3)(15,0)\put(7.5,
0){$\ldots$}\put(13,0){$\block$}\put(12.5, -1){\small$\phi(P'')$}
\end{picture}
\end{center}
\caption{Case 2.}\label{case2}
\end{figure}
\end{kst}
Conversely, given a matching $M$ on $[2n]$ which avoids both
patterns $12312$ and $121323$, we can get a Schr\"oder path $P$ of
semilength $n$ without peaks at level one. Suppose that
$M$ can be decomposed into smaller matchings
$\theta_1,\ldots,\theta_{k+1},\beta$ avoiding both patterns $12312$
and $121323$ as described in  Lemma~\ref{prop.1}.
If $k=0$ and $\theta_{1}=\emptyset$, then we have
$$\phi^{-1}(M)=H\phi^{-1}(\beta).$$
Otherwise, we get
$$\phi^{-1}(M)=U\phi^{-1}(\theta_1)UD\phi^{-1}(\theta_2)UD\ldots\phi^{-1}(\theta_{k})UD
\phi^{-1}(\theta_{k+1})D\phi^{-1}(\beta),$$
which is clearly a Schr\"oder path of semilength $n$ without peaks at level one.
Thus, we have constructed the desired bijection.

\begin{figure}[h,t]
\begin{center}
\begin{picture}(120,10)
\setlength{\unitlength}{4mm} \linethickness{0.4pt}
 \qbezier(15,0)(16.5,2)(18,0)\put(15,0){\circle*{0.2}}\put(14.5,-1){$1$}\put(18,0)
{\circle*{0.2}}\put(17.8,-1){$3$}\put(17,0){\circle*{0.2}}\put(16.5,-1){$2$}\put(27,0){\circle*{0.2}}
\put(27,-1){$12$}\qbezier(17,0)(22,3.5)(27,0)\put(19,0){\circle*{0.2}}\put(18.7,-1){$4$}\put(21,0)
{\circle*{0.2}}\put(21,-1){$6$}\qbezier(19,0)(20,2)(21,0)\put(20,0){\circle*{0.2}}\put(19.7,-1){$5$}
\put(24,0){\circle*{0.2}}\put(24,-1){$9$}\qbezier(20,0)(22,2.5)(24,0)\put(22,0){\circle*{0.2}}\put(21.8,-1){$7$}
\put(23,0){\circle*{0.2}}\put(23,-1){$8$}\qbezier(22,0)(22.5,
1.5)(23,0)\put(25,0){\circle*{0.2}}\put(24.6,-1){$10$}\put(26,0){\circle*{0.2}}\put(25.8,-1){$11$}
\qbezier(25,0)(25.5,1)(26,0)
\put(12.5,0){$\Longleftrightarrow$}
\put(0,0){\circle*{0.2}}\put(0,0){\line(1,1){1}}\put(1,1){\circle*{0.2}}\put(1,1){\line(1,1){1}}
\put(2,2){\circle*{0.2}}\put(2,2){\line(1,-1){1}}\put(3,1){\circle*{0.2}}\put(3,1){\line(1,-1){1}}
\put(4,0){\circle*{0.2}}\put(4,0){\line(1,1){1}}\put(5,1){\circle*{0.2}}\put(5,1){\line(1,1){1}}
\put(6,2){\circle*{0.2}}\put(6,2){\line(1,1){1}}\put(7,3){\circle*{0.2}}\put(7,3){\line(1,-1){1}}
\put(8,2){\line(1,-1){1}}\put(9,1){\line(1,0){2}}\put(11,1){\line(1,-1){1}}\put(12,0){\circle*{0.2}}
\put(8,2){\circle*{0.2}}\put(9,1){\circle*{0.2}}\put(11,1){\circle*{0.2}}
\put(2,-1.3){\small$UUDDUUUDDHD$}
\end{picture}
\end{center}
\caption{The bijection $\phi$.}\label{phi}
\end{figure}
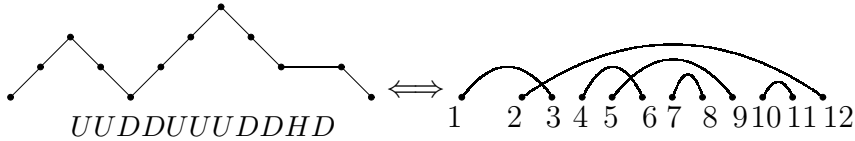

\begin{example} As illustrated in Figure \ref{phi},
  the Schr\"oder path $UUDDUUUDDHD$ corresponds to the matching
   $\{(1,3), (2,12), (4,6), (5,9), (7,8),
(10,11)\}$
\end{example}

In view of the bijection $\phi$, we see that a peak corresponds to
a crossing of the corresponded matching. Denote by
$\MM_{n,m}(12312, 121323)$ the set of the matchings in
$\MM_n(12312, 121323)$ with exactly $m$ crossings. We have the following formula.

\begin{theorem}
For $n,m\geq 0$, we have
  $$|\MM_{n,m}(12312, 121323)|={1\over n}{n\choose m}{2n-m\choose n+1}.$$
\end{theorem}

\pf It is well known that a Schr\"oder path of semilength $n$ can
be obtained from a Dyck path of semilength $n$ by turning some
peaks of the Dyck path into $H$ steps. A peak  is called a {\em
low} peak if it is at level one. Otherwise, it is called a {\em
high} peak. It has been shown  by Deutsch \cite{De} that  the
number of Dyck paths of semilength $n$  with exactly $k$ high
peaks is given by the Narayana number
\[ N(n,k)={1\over n}{n\choose k}{n\choose k+1}.\]
  Thus the number of Schr\"oder paths of semilength $n$
that contain exactly $m$ high peaks but no peaks at level one
equals
\begin{eqnarray*} \lefteqn{ \sum\limits_{k=0}^{n-1}{1\over
n}{n\choose k}{n\choose
k+1}{k\choose m}} \\[2mm]
& =  & \sum\limits_{k=0}^{n-1}{1\over n}{n\choose m}{n-m \choose
k-m}{n\choose k+1}\\[2mm]
& = & \sum\limits_{k=1}^{n}{1\over n}{n\choose m}{n-m \choose
n-k+1}{n\choose k}\\[2mm]
& = & {1\over n}{n\choose m}{2n-m \choose n+1}.
\end{eqnarray*}
This completes the proof.\qed

\section{\normalsize Matchings and Oscillating Tableaux}

In this section, we apply Stanley's bijection  between matchings and
oscillating tableaux to the set of $12312$-avoiding matchings \cite{st} to
obtain the corresponding restrictions on the oscillating tableaux for $12312$-avoiding
matchings. From the oscillating tableaux  we may
construct closed lattice walks and lattice paths that are counted by
the $3$-Catalan numbers.

It has been shown by Stanley~\cite{st} that oscillating tableaux
of length $2n$ are in one-to-one correspondence with matchings on
$[2n]$. We denote this bijection by $\rho$. This is stated as
follows. Given  an oscillating tableau $\emptyset=\lambda^{0},
\lambda^{1},\ldots,\lambda^{2n-1}, \lambda^{2n}=\emptyset$, we may
recursively define a sequence $(\pi_0, T_0), (\pi_1, T_1), \ldots,
$\break$(\pi_{2n}, T_{2n})$, where $\pi_i$ is
 a matching
and $T_i$ is a standard Young tableau (SYT). Let
$\pi_0$ be the empty matching and $T_0$ be the empty SYT.
The pair $(\pi_i, T_i)$  can be obtained from $(\pi_{i-1}, T_{i-1})$
 by the following procedure:
\begin{kst}
\item[(1)] If $\lambda^i\supset \lambda^{i-1}$, then
$\pi_i=\pi_{i+1}$ and $T_i$ is obtained from $T_{i-1}$ by adding
the entry $i$ in the square $\lambda^{i}\setminus \lambda^{i-1}$.

\item[(2)] If $\lambda^i\subset \lambda^{i-1}$, then let $T_i$ be
unique SYT of shape $\lambda^i$ such that $T_{i-1}$ is obtained
from $T_{i}$ by inserting some number $j$ by the RSK
(Robinson-Schensted-Knuth) algorithm. In this case, let
$\pi_i=\pi_{i-1}\cup (j,i)$.
\end{kst}
If the entry $i$ is added to $T_{i-1}$ to obtain $T_i$, then we
say that $i$ is added at step $i$. If $i$ is removed from
$T_{j-1}$ to obtain $T_j$, then we say that $i$ leaves at step
$j$. In this bijection,  $(i,j)$ is an edge of the corresponding
matching if and only if $i$ is added at step $i$ and leaves at
step $j$.
\begin{example}
For the oscillating tableau
     $$\emptyset, (1), (2), (2,1), (1,1), (1), \emptyset,$$
we get the corresponding sequence of $SYTs$ as follows:
$$
\begin{array}{lllllll}
\emptyset&\ \   1&\ \  12&\ \  12&\ \  1&\ \ 3&\ \ \
\emptyset,\\
&&&\ \  3&\ \  3&&
\\
\end{array}
$$
and the corresponding matching  with edges $\{(1,5), (2,4),
(3,6)\}$.
\end{example}

\begin{theorem}\label{coro.2}
There exists a bijection $\rho$ between the set of
$12312$-avoiding matchings on $[2n]$ and the set of oscillating
tableaux $T^\emptyset_{2n}$, in which each partition is of shape
$(k)$ or $(k,1)$ such that a partition $(k,1)$ is not followed immediately by
the partition $(k+1, 1)$.
\end{theorem}

\pf Let $M$ be a $12312$-avoiding matching. By definition, there
do not exist edges $(i_1, j_1)$, $(i_2, j_2)$ and $(i_3, j_3)$ in
$M$ such that $i_1<i_2<i_3<j_1<j_2$. Suppose that the
corresponding sequence of SYTs are
    $T_0, T_1, \ldots, T_{2n}$
and the  corresponding oscillating tableau  is
    $\lambda^0, \lambda^1, \lambda^2, \ldots, \lambda^{2n}$ under
    the Stanley's bijection.

If $\lambda^{p-1}$ is of shape $(k)$ for some $1\leq p\leq 2n$,
then it is possible that $\lambda^{p-1}\subset \lambda^{p}$. We
claim that all the entries in $T_{p-1}$ must leave the tableau in
deceasing order. It is clear that the entries in the first row of
an SYT are strictly decreasing from right to left. Suppose that
$i_2$ is right to $i_1$ in $T_{p-1}$ and $i_1$ leaves the tableau
before $i_2$. Assume that $i_1$ leaves at step $j_1$. That is to
say there exists an entry $h\geq p$ in the first row of  $T_{j_1}$
such that $T_{j_1-1}$ is obtained from $T_{j_1}$ by inserting
$i_1$. According to the RSK algorithm, the insertion of $i_1$
pushes $h$ up and $i_1$ takes the place of $h$.  Hence in
$T_{j_1}$, $h$ is left to $i_2$ in its first row, which
contradicts with the fact that $T_{j_1}$ is an SYT. So there do not
exist two crossing edges $(i_1, j_1)$ and $(i_2, j_2)$ such that
$i_1<i_2<p<j_1<j_2$. So $T_{p}$ can be obtained from $T_{p-1}$
by adding an entry $p$.

Now let us consider the case that $\lambda^{i-1}$ is
of shape $(k)$ and $\lambda^{i}$
is of shape $(k, 1)$. By the bijection $\rho$, $i$
is added  in the square $\lambda^i \setminus\lambda^{i-1}$
 to obtain $T_{i}$.
Suppose that $i$ is moved to the first row in $T_{j_1}$. Then
there exists a unique entry $j$ such that $T_{j_1-1}$ is obtained
from $T_{j_1}$ by row-inserting the entry $j$ by the RSK
algorithm. Hence $j<i$  and $j$ leaves before $i$. Suppose that
$i$ leaves the tableau at step $i_1$ with $j<i<j_1<i_1$, which
implies that $(j, j_1)$ and $(i, i_1)$ are two crossing edges of
the matching $M$. For any $i+1\leq p\leq j_1$,   we have
$\lambda^p\subset \lambda^{p-1}$. Otherwise, $T_p$ is obtained
from $T_{p-1}$ by inserting the entry $p$ in the square $\lambda^p
\setminus\lambda^{p-1}$.  It is obvious that $p$ is the initial
point of an edge of $M$. Let $(p, p_1)$ be an edge of $M$. Then
$(j, j_1)$, $(i, i_1)$ and $(p, p_1)$ are three edges of $M$ such
that $j<i<p<j_1<i_1$, which contradicts with the fact that $M$ is
a $12312$-avoiding matching. Furthermore, $\lambda^{j_1}$ is of
shape $(h)$ for some integer $h$. Thus we come to the assertion
that for the case when $\lambda^{p-1}$ is of shape $(k,1)$,  no
square is added to obtain $\lambda^{p}$ for any $1\leq p\leq 2n$.
\break This completes the proof.\qed

Given a matching $\pi\in \MM_n(12312)$, we may define a closed
lattice walk $\{\overrightarrow{v_i}=(x_i, y_i)\}_{i=0}^{2n}$
with $x_i\geq y_i$ in the $(x,y)$-plane from the origin to itself
by letting $x_i$ (resp. $y_i$) be the number of squares in the
first (resp. second) row of the partition $\lambda^i$ of the
corresponding oscillating tableau.
If $(x_{i+1}, y_{i+1})-(x_i, y_i)=(0,1)$, then by Theorem
\ref{coro.2} we see that the
size of  the next partition does not increase.
Thus we have the following corollary.

\begin{corollary} There is a one-to-one correspondence between
$12312$-avoiding matchings on $[2n]$  and
 closed lattice walks of length $2n$ in the
$(x,y)$-plane from the origin to itself consisting of the steps
$E=(1,0)$, $W=(-1, 0)$, $N=(0,1)$ and $S=(0,-1)$ such that a step
$N$ is  followed immediately by
 some consecutive $W$ steps and  one step $S$ and
 no step crosses the line $y=x$.
\end{corollary}

\begin{example}\label{ex.1}
The closed  lattice walk corresponding to the matching  $\{(1,5)$,
$(2,4)$, $(3, 6)\}$ is $EENWSW$.
\end{example}

 Denote by $L_n$ the set of  such closed lattice walks in the
$(x,y)$-plane as specified in the above corollary,
and denote by $P_n$ the set of   lattice
paths from $(0,0)$ to $(2n, n)$  consisting of  steps $E=(1,0)$
and $N=(0,1)$ and never  crossing the line $y={x\over 2}$. Hilton
and Pedersen \cite{hp} have shown that the cardinality of $P_n$ is
the $3$-Catalan number.
Now, let us describe a one-to-one correspondence between
$L_n$ and  $P_n$. Given a closed lattice walk $p\in L_n$,
we define a map $\tau$  by traversing the steps of $p$ along the
path and changing the steps of $p$ by the following rule:
$$\begin{array}{lll}
E &\rightarrow & EE,\\
W &\rightarrow & N,\\
N &\rightarrow & EN,\\
S &\rightarrow & E.
\end{array}$$
Denote by $|x|_p$ the number of $x$ steps in the path $p$.
Clearly, we  have that  $|E|_p=|W|_p$, $|N|_p=|S|_p$
and $ |E|_p+|N|_p=n$ since $p$ is a lattice path going from the
origin to itself with $2n$ steps. From the map $\tau$, we see that
 $|E|_{\tau(p)}=2|E|_p+|N|_p+|S|_p=2n$ and $ |N|_{\tau(p)}=|N|_p+|W|_p=n$.
 Hence $\tau(p)$ is a path  from $(0,0)$ to $(2n,n)$. We claim that
 $\tau(p)$  never crosses the line $y={x/ 2}$.   Otherwise,
  let $\tau(p)=p_1p_2\ldots p_{3n}$
 such that $p_k=N$ is the first step  going above the line $y={x/
 2}$. Let $p''=p_1p_2\ldots p_k$ and $\tau(p')=p''$. We get
$$|E|_{p''}=2|E|_{p'}+|N|_{p'}+|S|_{p'}<2|N|_{p''}=2|N|_{p'}+2|W|_{p'},$$
which implies that either $|E|_{p'}-|W|_{p'}<|N|_{p'}-|S|_{p'}$ or
$|E|_{p'}<|W|_{p'}$. This contradicts with the fact that  $p$ never goes
 above the line $y=x$, implying that $\tau(p)\in P_n$.

 Conversely,   given a lattice path $p\in P_n$, let $E_i$ be its
$i$-th $E$ step from left to right. If $E_{2k-1}$ and $E_{2k}$ are
consecutive steps in $p$, then $E_{2k-1}$ altogether with $E_{2k}$
corresponds to a $E$ step. Otherwise, $E_{2k}$ corresponds to
one $S$ step and $E_{2k-1}$ altogether with next $N$ step
corresponds to a $N$ step. For the remaining $N$ steps, each $N$
step corresponds to a $W$ step.

 Denote by $p'$ the resulted path. From the map we see that
 each $N$ step is followed by some $W$ steps and one
$S$ step in $p'$ and $|N|_{p'}=|S|_{p'}$. Moreover we have the
relations $|E|_{p}=2|E|_{p'}+|N|_{p'}+|S|_{p'}=2n$ and $
|N|_{p}=|N|_{p'}+|W|_{p'}=n$. It follows that $|E|_{p'}+|N|_{p'}=n$
and $|E|_{p'}=|W|_{p'}$, which implies that $p'$ is a path going
from the origin to itself with $2n$ steps. We claim that $p'$ is a
path never crossing the line $y=x$. Otherwise,  let $p'=p_1p_2\ldots
p_{2n}$,
    $p''=p_1p_2\ldots p_k$ and   $\tau^{-1}(p''')=p''$. From the map we have
$$|E|_{p'''}=2|E|_{p''}+|N|_{p''}+|S|_{p''}\geq 2|N|_{p'''}=2|N|_{p''}+2|W|_{p''},$$
and either $ |N|_{p''}=|S|_{p''} $ or $ |N|_{p''}-|S|_{p''}=1 $.
 Therefore, we obtain $|E|_{p''}-|W|_{p''}\geq |N|_{p''}-|S|_{p''}$,
 which implies that $p'$ is a path never crossing the line $y=x$.
 That is to say $p'\in L_n$. Up to now we have proved that
 the map $\tau$ is a bijection between  $L_n$ and  $P_n$, and we have the following
 conclusion.

\begin{theorem}
The map $\tau$ is a bijection between $L_n$ and $P_n$. Moreover,
we have
                  $$|L_n|=|P_n|=|\MM_n(12312)|={1\over 2n+1}{3n\choose n}.$$
\end{theorem}

\begin{example}
For $n=2$, we have
$$\begin{array}{cccc}
L_2:& EEWW &\ ENSW&\ EWEW\\
&\downarrow &\   \downarrow&\ \downarrow\\
 P_2:& EEEENN &\ EEENEN &\
EENEEN
\end{array}$$
\end{example}

\section{\normalsize Matchings and generating trees}

In this section we use the methodology of generating trees to deal
with other partial patterns  12132, 12123, 12321, 12231, 12213.
In fact, they are all Wilf-equivalent to 12312.
Given a matching $\pi$
on $[2n]$, a {\it position} $s$ of $\pi$ is meant to be the
position between the nodes $s$ and $s+1$ in the canonical sequential form
 if $1\leq s\leq 2n-1$,
and the {\it position} $2n$  is meant to be the position to the
right of the node $2n$. In the terminology of generating trees, a
position is called a {\em site}.

\begin{definition}
Let $\tau$ be a pattern  on $[k]$. The position $s$ of $\pi$ is an
active site if there exists a position $t$, $1\leq s\leq t\leq
2n$, such that inserting an edge starting at position $s$ and
ending at position $t$  gives a $\tau$-avoiding matching on
$[2n+2]$.  Otherwise,  the position $s$
 is said to be an inactive site.
\end{definition}

Given a partial pattern $\tau$, we use $T(\tau)$ to denote the generating tree
for the set of $\tau$-avoiding matchings on $[2n]$.

\begin{lemma}\label{lema1}
For any $\tau \in \{12312, 12132, 12123, 12321, 12231, 12213\}$,
the generating tree $T(\tau)$ is given by

\begin{equation}\label{eqtag1}
\left\{\begin{array}{l}\mbox{Root : }(0)\\[12pt] \mbox{Rule : }
(k)\rightsquigarrow(k+1)^1\,(k)^2\,(k-1)^3\ldots (0)^{k+2},
\end{array}\right.
\end{equation}
where the matching $\pi\in\MM_n(\tau)$ is labeled by $(k)$ such
that $k+2$ is the number of its active sites.
\end{lemma}

\pf We only consider the cases for $\tau=12312, \; 12123,\;
12321$. The other cases can be dealt with in the same manner.

The case $\tau=12312$: The matching $\pi$ on $[2]$ has two active
sites. This is consistent with the root label $(0)$.
Let $\pi$ be a $12312$-avoiding
matching on $[2n]$ labeled by $(k)$ with active sites
$i_1,i_2,\ldots,i_{k+2}$. Let $\pi'$ be a matching in
$\MM_{n+1}(12312)$ obtained from $\pi$ by inserting an edge from
position $i_s$ to position $i_t$ with $i_s\leq i_t$. Hence the
active sites of $\pi'$ are $i_t+1,i_t+2,
i_{t+1}+2,\ldots,i_{k+2}+2$. There are  $k+4-t$ such active sites.
So the children of the node $(k)$ are exactly the nodes
$(k')$ where $k'=k+2-t$. If $s$ ranges over  $1, 2, \ldots, k+2$
and $t$ ranges over $s, s+1, \ldots, k+2$, we get the desired rule
\ref{eqtag1}.

The case $\tau=12123$:  The proof is analogous to the previous
proof. The only difference lies in the following counting
argument. The active sites of $\pi'$ are
$i_s+1,i_{s+1}+1,\ldots,i_{t}+1$ when $i_s\neq i_t$ and
$i_s+1,i_{s}+2, i_{s+1}+2, \ldots,i_{k+2}+2$ when $i_s=i_t$. Thus,
the label of $\pi'$ is $(k')$ where $k'=t-s-1$ if $i_s\neq i_t$
and $k'=k-s+2$, otherwise. The rest of the argument is similar to
that in the previous case.

The case $\tau=12321$:   We only need to mention that the active
sites of $\pi'$ are $i_t+2, i_{t+1}+2, \ldots, i_{k+2}+2$ when
$t<{k+2}$ and $i_s+1, i_s+2, i_{s+1}+2 \ldots, i_{k+2}+2$ when
$t={k+2}$. It follows that the label of $\pi'$ is $(k')$ where
$k'=k-t+1$ when $t<k+2$ and $k'=k+2-s$ , otherwise.\qed

Applying Theorem \ref{theo1} and Lemma \ref{lema1}, we reach the
following  assertion.

\begin{theorem}\label{cbb}
For any $\tau \in \{12312, 12132, 12123, 12321, 12231, 12213\}$,
 we have $|M_{n}(\tau)|=C_{n,3}$.
\end{theorem}

 \vskip 2mm

\noindent{\bf Acknowledgments.} This work was done under the
auspices of the 973 Project on Mathematical Mechanization,
the National Science Foundation, the Ministry of
Education, and the Ministry of Science and Technology of China.


\end{document}